\documentclass[a4paper,oneside,12 pt]{amsart}
 
\usepackage[T1]{fontenc}
\usepackage[utf8]{inputenc}
\usepackage{amsthm}
\usepackage{amsmath}
\usepackage{amssymb}
\usepackage{mathrsfs}
\usepackage{graphicx}
\usepackage{enumitem} 
\usepackage{tikz}
\usepackage{color}
\usepackage{mathtools}
\usepackage[margin=1in]{geometry}

\newtheorem{definition}{Definition}[section]
\newtheorem{theorem}[definition]{Theorem}
\newtheorem{prop}[definition]{Proposition}
\newtheorem{remark}[definition]{Remark}

\newtheorem{lemma}[definition]{Lemma}
\newtheorem{cor}[definition]{Corollary}

\numberwithin{equation}{section}

\newcommand{\op}{\mathrm{op}}

\newcommand{\Cs}{\mathscr{C}}

\newcommand{\Zb}{\mathbb{Z}}

\newcommand{\Fb}{\mathbb{F}}

\newcommand{\Nb}{\mathbb{N}}
\newcommand{\Cb}{\mathbb{C}}

\newcommand{\Wb}{\mathbb{W}}

\newcommand{\D}{\mathcal{D}}

\newcommand{\C}{\mathcal{C}}
\newcommand{\M}{\mathcal{M}}
\newcommand{\N}{\mathcal{N}}

\newcommand{\U}{\mathcal{U}}
\newcommand{\R}{\mathcal{R}}

\newcommand\norm[1]{\left\| #1 \right\|}
\newcommand\md[1]{\left| #1 \right|}

\newcommand\p[1]{\left( #1 \right)}
\newcommand\set[1]{\left\lbrace #1 \right\rbrace}
\newcommand\floor[1]{\left\lfloor #1 \right\rfloor}

\newcommand{\e}{\varepsilon}

\begin{document}

\title{A metric characterization of freeness}

\author[]{L. Cadilhac}
\address{Laboratoire de Mathématiques d'Orsay, Univ. Paris-Sud, CNRS, Université Paris-Saclay, 91405 Orsay, France}
\email{leonard.cadilhac@u-psud.fr}

\author[]{B. Collins}
\address{Department of Mathematics, Kyoto University, Kitashirakawa Oiwake-cho, Sakyo-ku, 606-8502, Japan}
\email{collins@math.kyoto-u.ac.jp}

\date{}
\maketitle

\begin{abstract}
Let $\M$ be a finite von Neumann algebra and $u_1,\dots,u_N$ be unitaries in $\M$. We show that $u_1,\dots,u_N$ freely generate $L(\Fb_N)$ if and only if 
$$\norm{\sum_{i=1}^N u_i \otimes (u_i^{\op})^* + u_i^*\otimes u_i^{\op}}_{\M\overline{\otimes}\M^{\op}} = 2\sqrt{2N - 1}.$$
\end{abstract}

\section{Introduction}

The von Neumann conjecture formulated by Day in 1957 says that a group is not amenable if and only if it contains a non-amendable free group. It was first disproved by A. Ol'shanskii in 1980 \cite{Ols80}, and since then, the family of counterexamples has been expanded. 
A similar question can be asked at the level of von Neumann algebras: if a finite factor is not amenable, does it necessarily contain a free group factor? Little is known in that direction except for a breakthrough of Gaboriau and Lyons \cite{GL09}, who show that for certain wreath product groups $G$ (which may not contain $\Fb_2$), $L(\Fb_2) \subset L(G)$. Note that the authors of \cite{GL09} are, in fact, mainly interested in a version of the von Neumann conjecture for measure-preserving actions for which they provide a positive answer. 

A difficulty in tackling this problem is that there are no known abstract properties of $\M$ which would characterize the fact that $L(\Fb_2)$ embeds in $\M$. This remark first motivated us to write this note on a metric characterization of freeness. Although it is not clear that the results obtained (see Corollary \ref{main corollary}) can be used in the study of the von Neumann conjecture, we believe that they are of independent interest.
 
Indeed, they generalize at the operator level a well-known result of Kesten \cite{Kes59} who showed that given $g_1,\dots,g_N$ in a countable group $G$, the freeness of the $g_i$'s is characterized by the norm of the Markov operator associated to a random walk on $G$ supported by the $g_i$'s and their inverses. Let us denote by $\lambda: G \to L(G)$ the left-regular representation. From a von Neumann algebraic point of view, Kesten's result can be reformulated as follows: 
\begin{equation}\label{kesten}
g_1,\dots,g_N\ \text{are free in $G$}\ \Leftrightarrow 
\norm{\sum_{i=1}^N \lambda(g_i) + \lambda(g_i)^*}_{L(G)} = 2\sqrt{2N-1}.
\end{equation}

We extend this result by replacing the $\lambda(g_i)$'s with any finite family of unitary operators in a finite von Neumann algebra $\M$. The notion of freeness and Haar unitaries (unitaries with null moments) will be considered with respect to a fixed normal faithful tracial state $\tau$. We obtain the following:

\begin{theorem}\label{main}
Let $N \in \Nb$, $N>1$. Let $u_1,\dots,u_N$ be unitaries in $\M$. Then, the following assertions are equivalent:
\begin{enumerate}
\item the operators $u_1,\dots,u_N$ are free Haar unitaries,
\item the following equality holds:
$$\norm{\sum_{i=1}^N u_i \otimes (u_i^{\op})^* + u_i^*\otimes u_i^{\op}}_{\M\overline{\otimes}\M^{\op}} = 2\sqrt{2N - 1}.$$
\end{enumerate}
\end{theorem}

Note that $(1) \Rightarrow (2)$ is a consequence of $\eqref{kesten}$. Let us also mention that the inequality 
$$\norm{\sum_{i=1}^N u_i \otimes (u_i^{\op})^* + u_i^*\otimes u_i^{\op}}_{\M\overline{\otimes}\M^{\op}} \geq 2\sqrt{2N - 1},$$
is verified for any family of unitaries. 
This affirmation will become obvious when we reformulate the problem in terms of moments; it also appears in a paper of Pisier \cite{Pis97}. 
This leads to the following corollary:

\begin{cor}\label{main corollary}
Let $\M$ be a finite von Neumann algebra. Then the following are equivalent:\begin{enumerate}
\item $L(\Fb_2)$ embeds in $\M$,
\item $\inf_{u_1,u_2 \in \U(\M)} \norm{u_1 \otimes \bar{u_1}^* + u_1^*\otimes \bar{u_1} + u_2 \otimes \bar{u_2}^* + u_2^*\otimes \bar{u_2}}_{\M\otimes \M} = \sqrt{3},$
and this infimum is achieved. 
\end{enumerate}
\end{cor}

Let us make a few remarks in relation to this result. 
Firstly,  amenability can also be characterized by considering the same quantity. 
Indeed, in the spirit of Connes' \cite[Theorem 5.1]{Con76}, we know that a factor $\N$ is hyperfinite if and only if for every finite family of unitaries $u_1,\dots,u_n$ in $\N$, 
$$\norm{\sum_{i=1}^n u_i \otimes (u_i^*)^{\op}}_{\N \otimes^{\mathrm{min}} \N^{\op}} = n.$$
 In this sense, it is at the extreme opposite of freeness.
Secondly, it is worth pointing out that the problem considered in this manuscript complements 
results of \cite{Leh97}, and
also of \cite{Gri80,Coh82}, who consider generators of a group instead of general unitaries. 
Thirdly, in view of the above papers, it is natural to wonder what are the possible values of 
$$\norm{\sum_{i=1}^N u_i \otimes (u_i^{\op})^* + u_i^*\otimes u_i^{\op}}_{\M\overline{\otimes}\M^{\op}}$$
when the unitaries range on all possible choices in any tracial von Neumann algebra. 
This is clearly a subset of the interval $[2\sqrt{2N - 1}, 2N]$, and it can easily be seen that this is
the whole interval. Although many other approaches seem to be possible, let us
just outline one way to prove this assertion: take $N$ free unitary Brownian motions 
$\{t\mapsto u_i(t), i\in \{1,\ldots ,N\}\}$ as defined in \cite{Voi99}. Using explicit descriptions of the free unitary 
Brownian motion (see \cite{Bia97})
one can show that it is norm
continuous, converges to free Haar unitaries, and that this
convergence holds in norm for 
$t\mapsto \sum_{i=1}^N u_i(t) \otimes (u_i^{\op}(t))^* + u_i^*(t)\otimes u_i^{\op}(t)$, therefore
its norm is a continuous function taking value $2N$ at $t=0$ and tending to $2\sqrt{2N - 1}$. It follows
that the whole range $[2\sqrt{2N - 1}, 2N]$ is attained. 
Finally, it was brought to our attention by Franz Lehner that the traciality condition in Theorem \ref{main} is necessary, as illustrated by a counterexample appearing in his Ph.D. thesis \cite[p.51]{Leh97these}.

In section \ref{sec:combi}, 
we introduce our combinatorial 
approach to Theorem \ref{main}. Section \ref{sec:3}~contains the core technicalities: we use free group combinatorics to obtain a suitable lower bound on positive characters of $\Fb_N$, which allows us to conclude in Section \ref{sec:main}.

Acknowledgments: This work was initiated during the visit of LC to Kyoto University
early 2020, after preliminary discussions at MF Oberwolfach in 2018.
BC was supported by JSPS KAKENHI 17K18734 and 17H04823.
The authors are grateful to Mikael de la Salle, Cyril Houdayer, Eric Ricard, Adam Skalski, and Narutaka Ozawa
for inspiring comments and discussions.

\section{A combinatorial approach}\label{sec:combi}

\subsection{Reformulation of Theorem \ref{main}}

Let $N\in\Nb$. Let $\varphi$ be a positive definite function on the free group $\Fb_N$. We extend $\varphi$ linearly to $\Cb[\Fb_N]$ 
and keep the same notation, {\it i.e.} for any finitely supported family $(a_g)_{g\in\Fb_N}\in\Cb$,
$$\varphi\p{\sum_{g\in\Fb_N} a_g\cdot g} = \sum_{g\in\Fb_N} a_g\cdot \varphi(g).$$
Let $s_1,\dots,s_N \in \Fb_N$ be free generators of $\Fb_N$ and set 
$$a := \sum_{i=1}^N s_i + s_i^{-1} \in \Cb[\Fb_N].$$
We aim to prove the following:
\begin{theorem}\label{thm : main 2}
Assume that:
\begin{itemize}
\item $\varphi$ is constant on the conjugacy classes of $\Fb_N$ (it is a character), 
\item $\varphi(e) = 1$,
\item $\exists g\in\Fb_N, g\neq e, \varphi(g)\neq 0$,
\item $\forall g\in\Fb_N, \varphi(g)\geq 0$.
\end{itemize}
Then, 
$$\lim_{n\to\infty}\varphi(a^{2n})^{\frac1{2n}} > 2\sqrt{2N - 1}.$$
\end{theorem}
\begin{lemma}
Theorem \ref{thm : main 2} implies Theorem \ref{main}. 

\end{lemma}

\begin{proof}
	
\begin{sloppypar}
Consider the representation $\pi$ of $\Fb_N$ determined by \mbox{$\pi(s_i) = u_i \otimes (u_i^{\op})^*$} in $\U(\M \overline{\otimes} \M^{\op})$ for any $i \in \set{1,\dots,N}$. Define $\varphi := (\tau \otimes \tau^{\op}) \circ \pi$. Note that $\varphi$ is a positive character on $\Fb_N$ and that
\end{sloppypar}
\begin{align*}
\norm{\sum_{i=1}^N u_i \otimes (u_i^{\op})^* + u_i^*\otimes u_i^{\op}}_{\M\overline{\otimes}\M^{\op}} 
&= \lim_{n\to\infty} \norm{\sum_{i=1}^N u_i \otimes (u_i^{\op})^* + u_i^*\otimes u_i^{\op}}_{2n}\\
&= \lim_{n\to\infty}\varphi(a^{2n})^{\frac1{2n}}.
\end{align*}
Assume that the $u_i$'s are not free Haar unitaries. This means that there exists $g\in \Fb_N$ such that $g\neq e$ and $\varphi(g) \neq 0$. So $\varphi$ satisfies the conditions for Theorem \ref{thm : main 2} and hence 
$$\norm{\sum_{i=1}^N u_i \otimes (u_i^{\op})^* + u_i^*\otimes u_i^{\op}}_{\M\overline{\otimes}\M^{\op}}  > 2\sqrt{2N-1}.$$
\end{proof}

\subsection{Outline of the strategy of proof}

In this section, we give an overview of the strategy of proof. 
We start with the following remark:

\begin{remark}\label{rem:holder}
{\rm
To prove Theorem \ref{thm : main 2}, it suffices to show that 
$$\exists n, \varphi(a^{2n})^{\frac1{2n}} > 2\sqrt{2N-1},$$
by H\"older's inequality.
}
\end{remark}

More precisely, let $\tau$ be a normalized trace on $\M\overline{\otimes}\M^{\op}$, a simple
application of Stirling's formula shows that for any $N$ there exists $\varepsilon >0$ and $n_0$
such that for any $n\ge n_0$,
$$\tau (a^{2n}) > \varepsilon \cdot (2\sqrt{2N-1})^{2n}n^{-3/2}.$$

Looking into the assumptions of Theorem \ref{thm : main 2}, as soon as there is an element $g\in\Fb_N$, $g\neq e$ such that $\varphi(g)=\eta> 0$,
the traciality condition implies that all the conjugates of $g$ take the same trace value $\eta$, therefore, 
we can obtain a lower bound by observing the following:  given $N$, there exists $\varepsilon >0$ 
such that for $n$ large enough, the number of canonical paths 
on the free group on $N$ generators of length $n$, starting from the identity element and ending in a conjugacy class
is at least $\e n$ times the number of paths starting from the identity that end at a Cayley distance $\le 1$ from their starting point.

In turn, thanks to the positivity of the trace, this proves that 
there exists $\varepsilon >0$ and $n_0$
such that for any $n\ge n_0$,
$$\varphi (a^{2n}) > \varepsilon \cdot (2\sqrt{2N-1})^{2n}n^{-1/2}.$$

Here, we `gained' a factor $n$ in the subcritical part. 

This is not enough to conclude because $n^{-1/2}=o(1)$ and we would not be able to prove Theorem \ref{thm : main 2} without estimates on further conjugacy 
classes of elements other than $g$. 

Fortunately, there is a way to do so thanks to the fact that $\varphi$ is of positive type. Indeed, if
a finite collection of $g_i$ satisfy $\varphi (g_i)>\eta$ uniformly, then, assuming that the number of elements is large enough (for a threshold that depends on $\eta$), one can show that a positive proportion
with a ratio independent on the number of $g_i$'s)
 of the products $g_ig_j$ satisfies 
$\varphi (g_ig_j)>0$ uniformly as well. 
Thanks to this additional estimate and after picking the $g_i$'s a suitable amount of congujates of $g$, one can 
`gain' an additional factor $n$ in the subcritical part, and obtain the following
there exists $\varepsilon >0$ and $n_0$
such that for any $n\ge n_0$,
$$\varphi (a^{2n}) > \varepsilon \cdot (2\sqrt{2N-1})^{2n}n^{+1/2}.$$
Now, the exponent $+1/2$ is positive, so for $n$ large enough, the $L^n$ norm is bigger than $2\sqrt{2N-1}$,
which allows us to conclude thanks to Remark \ref{rem:holder}.

\subsection{Some notations}

We always assume $\Fb_N$ to be equipped with a set of distinguished generators $S = \set{s_1,\dots,s_N}$. 
Denote by $\Wb_N$ the set of words on the alphabet $S \cup S^{-1}$. For any $w \in \Wb_N$, denote by $l(w)$ its length. Words in $\Wb_N$ encode elements of $\Fb_N$ via the following map:

\begin{align*}
  g\colon \Wb_N & \to \Fb_N \\
  w_1\dots w_n & \mapsto w_1\cdots w_n.
\end{align*}
 
Recall that a word $w_1\dots w_n \in \Wb_N$ is said to be {\it reduced} if for any $i\leq n-1$, $w_i^{-1}\neq w_{i+1}$. The map $g$ is not injective but any $\gamma$ in $\Fb_N$ admits a unique reduced preimage. We will identify when convenient $\gamma$ with its writing as a reduced word $w$ in $\Wb_N$. In particular, set $l(\gamma):=l(w)$.

Define, for any $k \in \Nb$ and $w = w_1\dots w_k \in W_N$ of length $k$, the circular permutation of $w$ by
$$\sigma_1(w) = w_k w_1\dots w_{k-1}.$$
For any $t\in\Zb$, let $\sigma_t := \sigma_1^t.$ We say $w \in \Wb_N$ is {\it cyclically reduced} if $\sigma_t(w)$ is reduced for any $t$ 
(note that it is equivalent to $w$ and $\sigma_1(w)$ being reduced). If $w$ is a reduced word, $w$ can uniquely be written as 
$$w = u v u^{-1},$$
where $v$ is cyclically reduced $u$ is a reduced word. The element $v$ will be referred to as the {\it root} of $w$ and denoted by $r(w)$.

\section{Lower bound for positive characters of $\Fb_N$}\label{sec:3}

\subsection{Overview of the section}

Let $\varphi$ be as in Theorem \ref{thm : main 2} and $g_0\neq e$ such that $\varphi(g_0) = \alpha >0$. We will see soon (Corollary \ref{cor:evenlength}) that without loss of generality,
$g_0$ can be chosen to be of even length $2l_0$, $l_0 \in \Nb$. And since we might as well consider $r(g_0)$, it can be assumed to be cyclically reduced. This section is devoted to proving the following proposition:

\begin{prop}\label{prop:phik}
There exists a constant $C >0$ such that for any $k \geq l_0$,
$$\sum_{l(g) = 2k} \varphi(g) \geq Ck^2 (2N-1)^k.$$
\end{prop}

Note that the constant $C$ that we obtain above may be explicitly computed from $N$, $l_0$ and $\alpha$ but its precise value is inessential. 
Our goal is to exhibit a large enough family of words in $\Fb_N$ on which we can bound $\varphi$ from below. To do so, we combine 
two different operations:
\begin{itemize}
	\item The first operation consists of conjugating elements
	thanks to the tracial property. However, as explained in Section 2, by a simple calculation using the estimate given in Lemma \ref{lem:P2n0}, 
	it can be seen 	that considering only conjugates of $g_0$ is not  enough to obtain a sufficiently good lower 
	bound on $\varphi(a^{2n})$.
	\item The second operation is to multiply elements for which we already have a lower bound and apply Lemma \ref{lem : alpha beta}. 
	This is, however, a bit cumbersome because simplifications may occur, and estimating the length of the products obtained this way. 
	(to be able to use Lemma \ref{lem:P2n0}) requires some caution. 
\end{itemize}
The proof is divided into four different steps:
\begin{itemize}
\item Step 0 - we show that we can obtain new elements of positive trace by multiplication (Lemma \ref{lem : alpha beta}).
\item Step 1 - we fix an integer $i$ and apply Lemma \ref{lem : alpha beta} to the conjugates of $g_0$ of length $2i$ to obtain a set $\R_i$ new elements of quantified positive trace which have a length approximately $4i$ and are roots. The control on the length is essential to guarantee that our sets $\R_i$ do not overlap too much when $i$ varies. This will be essential during step $3$. 
\item Step 2 - we apply circular permutations to the elements of $\R_i$ to construct a bigger set $\R_i^{\sigma}$ of roots with positive traces.
\item Step 3 - we consider the sets $\Cs_i(k)$ of elements of length $2k$ with roots in $D_i$, show that they are essentially disjoint and big enough (in terms of the trace) to conclude. 
\end{itemize}

\subsection{Step $0$: locating more positivity through multiplication}

The following lemma allows obtaining new elements of positive trace through multiplication. 

\begin{lemma}\label{lem : alpha beta}
Let $k\in\Nb$ and $a>0$. Let $h_1,\dots,h_k$ be such that $\varphi(h_i) = \alpha$ for any $i\in \set{1,\dots,k}$. Then,
$$\sum_{1 \leq i,j \leq k} \varphi(h_i^{-1}h_j) \geq k^2\alpha^2.$$
\end{lemma}

\begin{proof}
Set $h_0 = e$ and consider the matrix
$$A = (\varphi(h_i^{-1}h_j))_{0\leq i,j\leq k}.$$
Since $\varphi$ is positive definite, $A$ is positive.
Let $\e >0$. Consider the vector $v = (1,-\e,\dots,-\e) \in \Cb^k$. Denote by $b$ the mean of  $(\varphi(h_i^{-1}h_j))_{0< i , j\leq k}$.
Note that
$$v^t A v = 1 - 2k\alpha\e + k^2b\e^2.$$
Set $\e = \frac1{k\alpha}$. By positivity of $A$, we have:
$$0 \leq 1 - 2 + \frac{b}{\alpha^2}.$$
Hence, $b\geq \alpha^2$.
\end{proof}

We can now justify that $g_0$ can be assumed to be of even length. 

\begin{cor}\label{cor:evenlength}
Let $\varphi$ be a character on $\Fb_N$ such that there exists $g_0 \neq e$ with $\varphi(g_0) > 0$. Then there exists $g_1$ of even length such that $\varphi(g_1) >0$. 
\end{cor}

\begin{proof}
Set $\alpha := \varphi(g_0)$. Assume that $g_0$ is of odd length, otherwise simply take $g_1 = g_0$. Choose $k$ large enough so that $\alpha^2 k^2 > k$ and $h_1,\dots,h_k$ distinct conjugates of $g_0$. Note that the $h_i$ are of odd length and hence the $h_i^{-1}h_j$ are of even length. 
\begin{align*}
\sum_{1\leq i\neq j \leq k} \varphi(h_i^{-1}h_j) &= \sum_{1\leq i,j \leq k} \varphi(h_i^{-1}h_j) - \sum_{1\leq i\leq k} \varphi(h_i^{-1}h_i) \\
&\geq \alpha^2k^2 - k \\
&> 0.
\end{align*}
This implies that there exists $i_1,j_1 \leq k$ such that $\varphi(h_{i_1}^{-1}h_{j_1}) > 0$. Set $g_1 = h_{i_1}^{-1}h_{j_1}$. 
\end{proof}

\subsection{Step $1$: products of conjugates of $g_0$}

The construction will depend on a parameter $k_0$. We could set a value for it right away, but it would not help with the clarity of the argument. We only need to know that it is large enough and only depends on $g_0$, $N$, and $\alpha$. 

Let $i \in \Nb$, $i > l_0$ and consider the set $\C_i$ of all conjugates of $g_0$ or $g_0^{-1}$ of length $2i$. Denote by $R_0$ the set of conjugates of $g_0$ or $g_0^{-1}$ of length $2l_0$. Note that since $g_0$ is assumed to be cyclically reduced:
\begin{equation}\label{eq:Ci}
	\md{\C_i} = \md{R_0}(2N-2)(2N-1)^{i - l_0 - 1} = c_0(2N-1)^i,
\end{equation}

where $c_0 = \md{R_0}(2N-2)(2N-1)^{-l_0-1}$. In the course of the proof, we will
apply circular permutations $\sigma_j$ to the words we obtain to generate new words. In order to guarantee at that future step that we indeed obtain new words, 
we remove from the beginning some pathological elements of $\C_i$.  For any $j\in\Nb$, define:
\begin{equation}\label{eq:defCij}
\C_i(j) := \set{g \in \C_i : \forall n,m \in \set{k_0+1,\dots,i-l_0} (\md{n-m} = 2j) \Rightarrow g_n = g_m},
\end{equation}
where $g_n$ denotes the $n$-th letter of $g$. 
Let $g \in \C_i(j)$, $g = w r_0 w^{-1}$ where $w$ is a reduced word, $r_0 \in R_0$ and no cancellation occur in this product. The condition above amounts to saying that $w$ is $2j$-periodic, except maybe for its first $k_0$ letters. Hence it is determined by the choice of its $k_0 + 2j$ first letters and consequently 
\begin{equation}\label{eq:Cij}
\md{C_i(j)} \leq c_0(2N-1)^{2j + k_0 + l_0}.
\end{equation}
Define 
\begin{equation}\label{eq:defCi'}
\C_i' := \C_i \backslash \bigcup_{j\leq \frac{i - l_0}{k_0}} \C_i(j).
\end{equation}
\begin{lemma}
For $i$ large enough and $k_0 \geq 4$, we have
\begin{equation}\label{eq:Ci'}
	\md{\C_i'} \geq 
	\dfrac{c_0}{2}(2N-1)^{i}.
\end{equation}
\end{lemma}

\begin{proof}
Since formula \eqref{eq:defCi'} only involves values of $j$ below $(i-l_0)/k_0$ and we assume $k_0 \geq 4$, we infer from \eqref{eq:Cij} that
$$
\md{\C_i(j)} \leq c_0(2N-1)^{k_0 + 2\frac{i-l_0}{k_0} + l_0} \leq K(2N-1)^{i/2},
$$
where $K = c_0(2N-1)^{k_0 + l_0/2}$. Hence, using \eqref{eq:Ci}, 
$$
\md{\C_i'} \geq \md{\C_i} - \dfrac{i-l_0}{k_0} K(2N-1)^{i/2} \geq \dfrac{c_0}{2}(2N-1)^i
$$
for $i$ large enough. 
\end{proof}

Define:
$$\R_i := \set{r(g) : g \in \C_i'\cdot\C_i', l(g) > 4i - 2k_0}.$$

\begin{remark}\label{rem:form}
{\rm
Let us keep in mind the form of the elements that we are dealing with. Let $h,h' \in \C_i'$. This means that there exists $r_0$ and $r_0'$ in $R_0$ and $u,u'$ in $\Fb_N$ such that $h = ur_0u^{-1}$ and $h' = u'r_0'(u')^{-1}$. Write $u = wv$ and $u' = wv'$ such that no cancellations occur in the product $v^{-1}v'$. Then $r(hh')$ belongs to $\R_i$ if and only if $l(w) < k_0$. In this case, $r(hh') = vr_0v^{-1}v'r_0'v'^{-1}$ and $l(r(hh')) = l(h) + l(h') - 4 l(w) \geq 4i - 4k_0$. 
}
\end{remark}

\begin{lemma}\label{lem:Ci2}
For $i$ and $k_0$ large enough, the following estimate holds:
$$ \sum_{g\in\R_i} \varphi(g) \geq \dfrac{c_0}{8} \alpha^2	(2N - 1)^{2i - k_0}.$$
\end{lemma}

\begin{proof}
For any $g \in \C_i'$, there are at most $c_0(2N-1)^{i-k_0}$ elements $g'$ of $\C_i'$ for which $l(gg') \leq 4i - 2k_0$. 
Indeed, the first $k_0$ letters of $g'$ must coincide with the first $k_0$ letters of $g$ in order for at least $k_0$ simplifications to occur. This means that:
	
$$\sum_{g, g' \in  \C_i'} \varphi(gg') \leq c_0(2N-1)^{i-k_0}\md{\C_i'} 
+ \sum_{\substack{g, g' \in  \C_i'\\ l(gg') > 4i - 2k_0}} \varphi(gg').$$
Then, by Lemma \ref{lem : alpha beta} 
(note that considering $gg'$ or $g^{-1}g'$ does not modify the following sum since $\C_i' = (\C_i')^{-1}$):
\begin{align*}
	\sum_{\substack{g, g' \in  \C_i'\\ l(gg') > 4i - 2k_0}} \varphi(gg') &\geq
	\md{\C_i'}^2\alpha^2 - c_0(2N-1)^{i-k_0}\md{\C_i'}\\
	&\geq \dfrac{c_0^{2}\alpha^2}{4}(2N-1)^{2i} - \dfrac{c_0^2}{2}(2N-1)^{2i-k_0} \\
	&\geq \dfrac{c_0^{2}\alpha^2}{8}(2N-1)^{2i},
\end{align*}
where we assume in the last line that $k_0$ is large enough so that,
$$(2N-1)^{-k_0} \leq \dfrac{\alpha^2}{4}.$$
Finally, note that at most $(2N-1)^{k_0}$ pairs $h,h' \in \C_i'$ can give rise to the same element $r(hh') \in \C_i^{(2)}.$ This means that 
$$ (2N-1)^{k_0} \sum_{g\in\R_i} \varphi(g) \geq \sum_{\substack{g , g' \in  \C_i'\\ l(gg') > 4i - 2k_0}} \varphi(gg'),$$
which is the desired estimate.
\end{proof}

\subsection{Step 2: circular permutations}

To guarantee that we construct new elements by circular permutations, we rely on the following observation.

\begin{lemma}\label{lem:circular}
Let $g,g' \in \Fb_N$ and $t\in \Nb$. Assume that $g$ and $g'$ can be written as follows
\begin{center}
$g = uau^{-1}b$ and $g' = u'a'(u')^{-1}b'$ with $l(a) = l(a'), l(u) = l(u'), l(b) = l(b')$,
\end{center}
and no cancellation in the products. Then if $g = \sigma_t(g')$, $u$ is $2t$-periodic. 
\end{lemma}

\begin{proof}
Denote by $l$ the length of $u$. We may assume that $2t < l$ otherwise the statement is void. For $i \in \set{t+1,\dots,l}$, note that the $i$-th letter of $g$ is $u_i$ and the $i$-th letter of $\sigma_t(g')$ is $u'_{i-t}$. So 
\begin{center}
$u_i = u'_{i-t}$ for $i\in \set{t+1,\dots,l}$.
\end{center}
Now we look at the occurrences of $u^{-1}$ and $(u')^{-1}$ in $g$ and $\sigma_t(g')$ {\it i.e.} the $i$-th letters of $g$ and $\sigma_t(g')$ for $i \in \set{l + l(a) + t+1,\dots,l+l(a)+l}$. We obtain 
\begin{center}
$u_i = u'_{i+t}$ for $i \in \set{1,\dots,l-t}.$
\end{center}
Combining the two equalities, we get that $u_i = u_{2t+i}$ for $i \in \set{1,\dots,l-2t}.$  
\end{proof}

Define:
$$\R_i^\sigma := \set{\sigma_j(g): 0 \leq j \leq \dfrac{i-l_0}{k_0}, g\in \R_i }.$$

\begin{lemma}\label{lem:Di}
	For $i$ and $k_0$ large enough, the following estimate holds:
	$$ \sum_{g\in\R_i^\sigma} \varphi(g) \geq \floor{\dfrac{i - l_0}{k_0} }\dfrac{c_0}{8} \alpha^2	(2N - 1)^{2i - k_0}.$$
\end{lemma}

\begin{proof}
Given the estimate obtained in Lemma \ref{lem:Ci2}, it suffices to prove that for any $0\leq j,j' \leq (i-l_0)/k_0$ and $g,g' \in \R_i$,
$$\sigma_j(g) = \sigma_{j'}(g') \Rightarrow (g = g'\ \text{and} \ j = j').$$
Let us assume that $\sigma_j(g) = \sigma_{j'}(g')$. If $j = j'$ then immediately $g = g'$ and we get the expected conclusion. So assume by contradiction that $j\neq j'$. Set $t = j'-j$, we have 
$$
g = \sigma_{-j}\circ \sigma_{j}(g) = \sigma_{-j}\circ\sigma_{j'}(g') = \sigma_t(g').
$$
According to Remark \ref{rem:form}, $g$ and $g'$ can be written as follows:
$$g = ur_0u^{-1}vs_0v^{-1} \quad g' = u'r_0'(u')^{-1}v's_0'(v')^{-1}$$
where $r_0,s_0,r_0',s_0'$ belong to $R_0$ and no cancellations occur in the products above. Still by Remark \ref{rem:form}, $u$ and $v$ have the same length (same for $u'$ and $v'$) and since $l(g) = l(g')$, we get that $u,u',v,v'$ are words of the same length $l$. By Lemma \ref{lem:circular}, $u$ is $2t$-periodic. 
Since $g$ is an element of $\R_i$, $g$ comes from a product $hh'$, $h,h' \in \C_i'$. Write $h = wv_0w^{-1}$. We have $w = w_1\dots w_{i-l_0-l} u$. By construction of $\C_i'$, $w$ cannot be $2t$-periodic starting from its $k_0$-th letter, which is a contradiction (see \eqref{eq:defCi'}, \eqref{eq:defCij}, and recall that by assumption $t \leq \dfrac{i-l_0}{k_0}$).
\end{proof}

\subsection{Step 3: conjugation}

For any $k \geq 2i$, define: 
$$\Cs_i(k) := \set{g \in \Fb_N : l(g) = 2k, r(g) \in \R_i^{\sigma}}.$$
\begin{lemma}\label{lem:Dik}
Let $i,i'$ large enough and assume that $\md{i-i'} \geq k_0$. Let $k \geq 2i, 2i'$. Then,
\begin{itemize}
\item $\Cs_i(k)$ and $\Cs_{i'}(k)$ are disjoint,
\item There exists $C'>0$ independent of $k$ or $i$ such that:
$$\sum_{g\in \Cs_i(k)} \varphi(g) \geq C' i (2N-1)^k.$$
\end{itemize}

\begin{proof}
For the first point, recall that by Remark \ref{rem:form}, the elements of $\R_i$ have length between $4i$ and $4i - 4k_0 +4 $. Since those elements are roots, applying circular permutations to them does not change their length. Hence, the elements of $\R_i^{\sigma}$ also have length between $4i$ and $4i - 4k_0 +4 $. So if $\md{i-i'} \geq k_0$ then $\R_i^{\sigma}$ and $\R_{i'}^{\sigma}$ are disjoint so $\Cs_i(k)$ and $\Cs_{i'}(k)$ are disjoint.

Let $g$ be a cyclically reduced element of length $4i$ in $\Fb_N$. Note for $k> 2i$, there are $(2N-2)(2N-1)^{k-2i-1}$ elements of $\Fb_N$ of length $2k$ and root $g$. This means that any element of $\D_i$ is the root of at least $(2N-2)(2N-1)^{k-i-1}$ elements in $\Cs_i(k)$. Hence,
$$\sum_{g\in \Cs_i(k)} \varphi(g) \geq (2N-2)(2N-1)^{k-2i-1} \sum_{g\in \R_i^{\sigma}} \varphi(g).$$
By Lemma \ref{lem:Di}, we obtain the expected estimate.
\end{proof}
\end{lemma}

\begin{proof}[Proof of Proposition \ref{prop:phik}]
Note that for any $k \geq l_0$, $\sum_{l(g) = 2k} \varphi(g) >0$ so it suffices to prove Proposition \ref{prop:phik} for large values of $k$. Fix a rank $i_0$ such that for $i,i'\geq i_0$, Lemma \ref{lem:Dik} applies. Let $k \in \Nb$, $k \geq 2(i_0 + k_0)$ (so that the sums below are non empty). By Lemma \ref{lem:Dik}, we have:
\begin{align*}
\sum_{\md{g} = 2k} \varphi(g) 
&\geq \sum_{\substack{i_0 \leq i \leq k/2 \\ k_0 | i}} \sum_{g\in \Cs_i(k)} \varphi(g) \\
&\geq \sum_{\substack{i_0 \leq i \leq k/2 \\ k_0 | i}} C'i(2N-1)^k \\
&\geq Ck^2 (2N-1)^{2k},
\end{align*}
where $C$ is a small enough constant independent of $k$. 
\end{proof}

\section{Proof of the main theorem}\label{sec:main}

\begin{lemma}\label{lem:P2n0}
Let $n\in\Nb$ and $k<n$. Let $g\in\Fb_N$ such that $l(g) = 2k$. Then:
$$\md{\set{w\in\Wb_N
 : l(w) = 2n, g(w) = g}} 
 =
 \p{\binom{2n}{n-k} - \binom{2n}{n-k-1}}(2N-1)^{n-k} =: N_{n,k}.$$
\end{lemma}

\begin{proof}
We interpret words as paths on the Cayley graph of $\Fb_N$. Note that the quantity we consider only depends on $l(g)$. A way to generate paths going from $e$ to {\it any} element of length $2k$ is first to choose at which times the path is going to go away from $e$ and at which times the path is going to come back towards $e$. Since the path is going 
to an element of length $2k$, it has to go $n+k$ times away from $e$ and $n-k$ back to $e$. The number of possible choices there, for a path of 
length $2n$, is given by the Catalan triangle 
$$C_{n,k} = \binom{2n}{n-k} - \binom{2n}{n-k-1}.$$
Moreover, when we chose to go away from $e$, there are {\it at least} $2N-1$ possible directions, and $2N$ possible directions for the first time, 
thus obtaining
$$C_{n,k}2N(2N-1)^{n+k-1}$$
paths. Finally, since, for now, all we have fixed is the length of the target of the path and not a particular point, we have to divide this result by the 
number of elements of length $2k$ in $\Fb_N$ {\it i.e.} $2N(2N-1)^{2k-1}$, to get the desired estimate.
\end{proof}

We are now ready to prove our main theorem. 

\begin{proof}[Proof of Theorem \ref{thm : main 2}]
Let $n\in\Nb$. First remark that by Lemma \ref{lem:P2n0} and Proposition \ref{prop:phik},
\begin{align*}
\varphi(a^{2n}) &= \sum_{w\in\Wb_N} \varphi(g(w)) \\
&\geq \sum_{k\leq n} C_{n,k}(2N-1)^{n-k} \sum_{\md{g} = 2k} \varphi(g) \\
&\geq C\sum_{l_0 \leq k\leq n} C_{n,k} k^2 (2N-1)^n.
\end{align*}
Note that given the expression of $C_{n,k}$,
$$ \sum_{l_0 \leq k \leq n}k^2 C_{n,k} \geq
 \sup_{l_0 \leq k \leq n} k^2 \sum_{k\leq i \leq n} C_{n,k} = \sup_{l_0 \leq k \leq n} k^2 \binom{2n}{n-k}.$$
Now chose $k = n^{1/3}$. For $n$ large,
\begin{align*}
\binom{2n}{n-k} &\sim \dfrac{2^{2n}}{\sqrt{\pi n}}\p{\dfrac{n}{n+k}}^{n+k}\p{\dfrac{n}{n-k}}^{n-k} \\
&\sim \dfrac{2^{2n}}{\sqrt{\pi n}}e^{-k + o(1)}e^{k+o(1)} \\
&\sim \dfrac{2^{2n}}{\sqrt{\pi n}}.
\end{align*}
Hence, for $n$ large enough,
$$\varphi(a^{2n}) \geq \dfrac{C}{\sqrt{\pi}} n^{1/6}2^{2n}(2N-1)^n > 2^{2n}(2N-1)^n.$$
Which concludes the proof by Remark \ref{rem:holder}.
\end{proof}

\begin{proof}[Proof of Corollary \ref{main corollary}]
Let $\tau$ be a faithful trace on $\M$ and consider $L^2(\M,\tau)$. It follows from 
Theorem \ref{main} that all reduced words in $u_i$ and their inverses form 
an orthonormal family. Since the GNS representation of the von Neumann subalgebra
generated by $u_1,u_2$ is faithful; it allows to conclude that is is isomorphic to $L(\Fb_2)$.
\end{proof}

\bibliography{Bibli}
\bibliographystyle{plain}

\end{document}